\documentclass[a4paper,12pt]{article}
\usepackage[english]{babel}
\usepackage[utf8]{inputenc}
\usepackage{amsmath}
\usepackage{amssymb}
\usepackage{amsfonts}
\usepackage{graphicx}
\usepackage{color}
\newtheorem{Theorem}{Theorem}
\newtheorem{Remark}{Remark}
\begin{document}

\title{On symmetric wedge mode of an elastic solid}

\author{Alexander Nazarov$^{ab}$\footnote{Corresponding author; e-mail: al.il.nazarov@gmail.com}, Sergey Nazarov$^{ac}$, and German Zavorokhin$^b$}
\date{}
\maketitle

\vspace{-6mm} 

\begin{center}
{$^a$\it Department of Mathematics and Mechanics,
St.Petersburg State University,\\
Universitetskii prospect, 28, 198504,
St.Petersburg, Russia \\
%e-mail: al.il.nazarov@gmail.com\\
$^b$St.Petersburg Department of the Steklov Mathematical Institute,\\
Fontanka, 27, 191023, St.Petersburg, Russia \\
%e-mail: zavorokhin@pdmi.ras.ru\\
%e-mail: srgnazarov@yahoo.co.uk
$^c$Institute for Problems in Mechanical Engineering of RAS,\\
Bol'shoy Prospekt V.O., 61, St. Petersburg, 199178, Russia\\
}

\end{center}

\date{}
\maketitle

\def\e{\varepsilon}

\bigskip \noindent {\bf Abstract.} 
The existence of a symmetric mode in an elastic solid wedge for all admissible values of the Poisson ratio and arbitrary interior angles close to $\pi$ has been proven.

\bigskip

\hfill {\it To Professor V.M. Babich on the occasion on his 90th anniversary, with admiration}

\section{Introduction. Statement of the problem}

Along with bulk and surface waves (a fair example of the latter is the Rayleigh wave \cite{LR} in half-space), wedge waves comprise a fundamental type of oscillations of solids and are intensively studied in geophysics, machine building, civil engineering etc.

First theoretical results on waves propagating along the edge of a wedge were obtained by numerical simulations (see, e.g., \cite{La}, \cite{La1}, \cite{La2}, \cite{MMC}, \cite{TR}). Then these waves were studied analytically at the physical level of rigor by many authors, mostly for small interior angles (slender wedge, see, e.g., \cite{Mc}, \cite{Mo}, \cite{Kr1}, \cite{Kr}, \cite{Kr2}, \cite{P}) and interior angles slightly less then  $\pi$ (\cite{Mozh1}, \cite{Mozh}, \cite{P}, \cite{Sh}, see also a survey in \cite[Ch. 10, Sec. 2]{BGK}). 

The first rigorous proof of existence of the wedge wave was obtained in the pioneering paper \cite{K} by a variational method for interior angles less then $\frac \pi 2$. Then the idea of \cite{K} was developed in \cite{ZN} 
and \cite{Pu} where the range of aperture angles was enlarged. Further applications of this idea can be found in \cite{Ba} and \cite{PLNM}.

In this paper, we prove the existence of a symmetric wedge mode in an elastic deformable wedge for all admissible values of the Poisson ratio $\sigma\in(-1,\frac 12)$ and interior angles close to $\pi$ and derive an 
asymptotic formula for corresponding eigenvalue. We note that neither explicit coefficient in the asymptotic formula for the velocity of the wave propagating along elastic wedges with interior angles close to $\pi$, 
nor a validation of the formal computations were provided in above-mentioned papers, and our approach is novel and differs from applied previously.

\smallskip

\begin{figure}[h!]
	\begin{center}
			\includegraphics[width=0.7\linewidth]{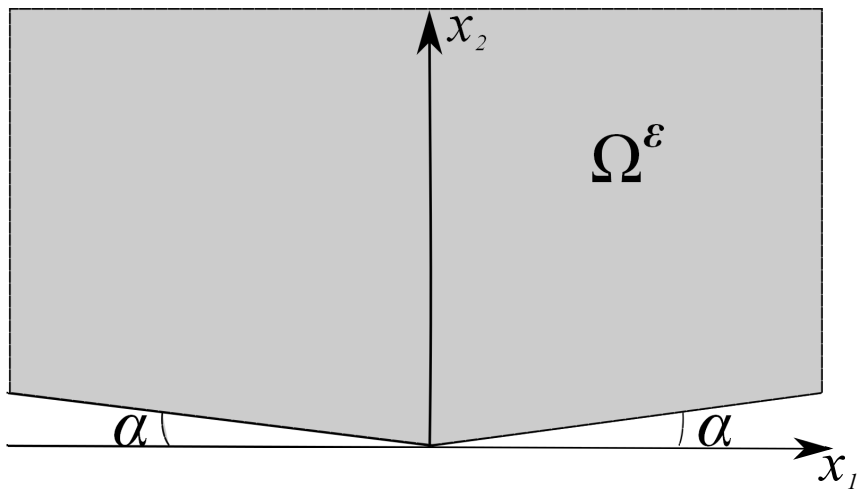}
			\caption{Wedge cross section}  
			\label{Fig1}  
	\end{center} 
\end{figure}

Let
$$
  \Omega^\e=\{(x_1,x_2): r>0, \ \phi\in (\alpha,\pi - \alpha)\}
$$
be an angular domain of the plane $\mathbb{R}^2$. Here $(r,\phi)$ stand for polar coordinates of a point $(x_1, x_2)$, and $\e=\tan(\alpha)$ is a small positive parameter. We consider $\Omega^\e$ as the cross section of 
an isotropic homogeneous elastic wedge $K^\e=\Omega^\e \times \mathbb{R} \ni (x_1,x_2,x_3)$; see Fig. \ref{Fig1}. 

We seek for a solution ${\cal U}^\e=({\cal U}_1^\e,{\cal U}^\e_2,{\cal U}_3^\e)$
of the elasticity system in $K^\e$ subject to the traction-free boundary conditions 
on $\partial K^\e$ (see \cite[Ch. 3]{LL}) in the form
\begin{equation*}\label{mode}
{\cal U}^\e(x_1,x_2,x_3,t)=
{\mathbb U}^\e(x_1,x_2,x_3)e^{-i \omega t},
\end{equation*}
where $\omega=\omega(\e)>0$ is a given oscillation frequency. Then for the vector function of displacements $ {\mathbb U}^\e = ({\mathbb U}_ {1}^\e, {\mathbb U}_ {2}^\e, {\mathbb U}_ {3}^\e) $, we obtain the 
boundary value problem 
\begin{equation}\label{2dim1}
 {\cal L}(\partial_1,\partial_2,\partial_3)
 {\mathbb U}^\e  =  \rho \omega^2 {\mathbb U}^\e \quad\mbox{in}\quad K^\e,\qquad
%\end{equation}
%\begin{equation}\label{2dim2}
{\cal N}(\partial_1,\partial_2,\partial_3) {\mathbb U}^\e  =   0\quad\mbox{on}\quad\partial K^\e,
\end{equation}
where $\rho>0$ is the constant density of the material (%in the sequel, 
without loss of generality, we set $\rho=1$), and $\partial_m=\frac{\partial}{\partial x_{m}}$.

Differential expressions in (\ref{2dim1}) 
are defined component-wise by the following relations with $n=1,2,3$:
\begin{equation}\label{sigma}
\big({\cal L}(\partial_1,\partial_2,\partial_3) {\mathbb U}^\e\big)_n=-\,\partial_m\sigma_{nm}
[{\mathbb U}^\e],
\qquad
\big({\cal N}(\partial_1,\partial_2,\partial_3) {\mathbb U}^\e\big)_{n}=\sigma_{nm}[{\mathbb U}^\e]\nu_{m},%\quad n=1,2,3,
\end{equation}
where $\nu=(\nu_{1},\nu_{2},\nu_{3})$ is the outward normal unit vector on $\partial K^\e$ (evidently $\nu_3=0$).
Here and throughout the paper we use the convention of summation from $1$ to $3$ with respect to repeated indices. 

In (\ref{sigma}) $\sigma_{nm}[{\mathbb U}^\e]$ stand for the stress tensor components
\begin{equation*}\label{formula}
\sigma_{nm}[{\mathbb U}^\e]=\lambda\,\partial_\ell {\mathbb U}_{\ell}^\e\delta_{n,m} +\mu \big(\partial_m{\mathbb U}_{n}^\e +\partial_n {\mathbb U}_{m}^\e\big),
\end{equation*}
where $\delta_{n,m}$ is the Kronecker symbol and the Lam\'e constants $\lambda$ and $\mu$ satisfy the inequalities $\lambda+\frac{2}{3}\mu>0$ and $\mu>0$. 

\begin{Remark}\label{rotations}
Notice that the equation in (\ref{2dim1}) can be rewritten in the form
\begin{equation}\label{div}
-\mu\Delta {\mathbb U}^\e -(\lambda+\mu){\rm grad}\ {\rm div}\,{\mathbb U}^\e=\omega^{2}{\mathbb U}^\e, 
\end{equation}
while the boundary condition in (\ref{2dim1}) is just the natural boundary condition for (\ref{div}). Therefore, the problem (\ref{2dim1}) is invariant under rotations.
\end{Remark}

We search for a solution to the problem (\ref{2dim1}) %, (\ref{2dim2})
in the form 
\begin{equation*}
{\mathbb U}^\e(x_1,x_2,x_3)=u^\e(x_1,x_2)e^{i kx_3},
\end{equation*}
corresponding to an elastic wave propagating along the edge of the wedge $K^\e$ while $k>0$ is 
the wave number. Then the problem reduces to
\begin{equation}\label{2dim11}
{\cal L}(\partial_{1},\partial_{2},ik)u^\e  =  \omega^2 u^\e\quad\mbox{in}\quad\Omega^\e,\qquad
%\end{equation}
%\begin{equation}\label{2dim22}
{\cal N}(\partial_{1},\partial_{2},ik)u^\e  =   0\quad\mbox{on}\quad\partial\Omega^\e.
\end{equation}
From now on we regard that in all formulae $\partial_{3}$ should be replaced by $ik$.

Since the wedge waves are localized in the vicinity of the edge, the vector function $ u^\e = (u_ {1}^\e, u_ {2}^\e, u_ {3}^\e) $ should decay at infinity.

Due to material symmetry with respect to the plane $x_1=0$, the operator of the problem (\ref{2dim11}) is reduced by the decomposition of the Sobolev space 
$\big(H^{2}(\Omega^\e)\big)^3$ 
into the subspaces of antisymmetric and symmetric displacements, which makes it possible to consider the antisymmetric (flexural) modes
\begin{equation*}\label{antisymm}
%\left.\begin{array}{ccc}
u^{\e}_{1}(-x_1, x_2)=u^{\e}_{1}(x_1,x_2), \quad %\\
u^{\e}_{2}(-x_1, x_2)=-u^{\e}_{2}(x_1,x_2), \quad %\\
u^{\e}_{3}(-x_1, x_2)=-u^{\e}_{3}(x_1,x_2), %\\
%\end{array} \right.
\end{equation*}
and the symmetric modes
\begin{equation}\label{symm}
%\left.\begin{array}{ccc}
u^{\e}_{1}(-x_1, x_2)=-u^{\e}_{1}(x_1,x_2), \quad % \\
u^{\e}_{2}(-x_1, x_2)=u^{\e}_{2}(x_1,x_2), \quad % \\
u^{\e}_{3}(-x_1, x_2)=u^{\e}_{3}(x_1,x_2). %\\
%\end{array} \right.
\end{equation}
 
Numerical simulations in \cite{MMC}, \cite{TR} predict the existence of a symmetric mode (\ref{symm}) for obtuse wedges. The existence of such mode was proved rigorously in \cite{ZN} and \cite{Pu} for certain 
values of the Poisson ratio $\sigma=\frac{\lambda}{2(\lambda+\mu)}$ and some range of wedge interior angles that are far from $\pi$.

Now we formulate our main result.

\begin{Theorem}\label{main}
For any $\sigma\in(-1,\frac 12)$, one finds $\e^0$
such that for any $0<\e<\e^0$ there exists a symmetric mode $u^\e \in \big(H^2(\Omega^{\e})\big)^3$ decaying exponentially as $|x|\to\infty$ and 
solving the problem (\ref{2dim11}) with $\omega^2=\widehat\omega^2(\e)<{\bf c}_R^2k^2$. Corresponding wedge wave propagates along the edge with the velocity ${\bf c}_w$ which has the following asymptotics as $\e\to0$:
\begin{equation}\label{asymp-velo}
{\bf c}_w^2={\bf c}_R^2(1-\e^2\vartheta+O(\e^{\frac 52})),
\end{equation}
where $\vartheta>0$ is an explicit coefficient depending on $\sigma$ only.
\end{Theorem}

A short announcement of our results was published in \cite{ZNN}. Unfortunately, the text in \cite{ZNN} contains essential errors and misprints.\medskip

Recall that rigorous results obtained in previous papers were established by variational methods. However, a similar approach is not applicable for interior angles close to $\pi$, see \cite{SN2010}. In this way, we use an asymptotic method, 
namely, a variant of the method of matched asymptotic expansions (see \cite{vD}, \cite{Il}, \cite[Ch. 2]{MNP} and other monographs) adapted earlier only for scalar diffraction problems in \cite{SN2010}, \cite{SN2011} and other papers. 
\medskip 

The paper is organized as follows. First of all, in Section \ref{trans} we use the isotropy of the elastic material to transform an angular domain (Fig. \ref{Fig1}) to a half-plane with a narrow wedge-shaped notch (Fig. \ref{Fig3}) which disappears upon formal passage to the limit $\e\to+0$ (cf. \cite{SN2011} for a scalar problem). Here we also give the operator formulation of the problem. 

In Subsection \ref{Ss:pencil} we apply the Fourier--Laplace transform with the dual variable from a thin strip containing the imaginary axis and study the spectrum of the resulting quadratic pencil. 
Then, within the framework of the method of matched asymptotic expansions, well-known asymptotic procedures are applied. Namely, the ``outer'' expansion (Section~\ref{outer}) involves the asymptotics of the eigenvalues and eigenvectors under perturbation of operator pencils, while the ``inner'' expansion (Section \ref{inner}) reproduces the constructions of solutions in domains with regularly perturbed boundary, as well as elements of the Kondrat'ev theory \cite{Ko}. We emphasize that this theory is applied in a non-standard situation: a half-plane is interpreted as a strip of infinite width, see Subsection \ref{Kondr}. That is why preliminary work of Subsection \ref{Ss:pencil} was required. 

In Section \ref{Justify} we construct the approximate eigenvector of the problem (\ref{2dim11}) gluing the asymptotic expansions obtained in Sections \ref{inner} and \ref{outer}. Using elements of the theory of spectral measure, we justify the asymptotic formula (\ref{asymp-velo}) and prove Theorem \ref{main}.
\medskip

We recall the following standard notation:

${\bf c}_t=\sqrt{\mu}$ and ${\bf c}_l=\sqrt{\lambda+2\mu}$ are the velocities of the transverse and longitudinal bulk waves, respectively;

${\bf c}_R<{\bf c}_t$ is  the  velocity  of  the  Rayleigh  wave;

$\varkappa_{t}^2=1-{{\bf c}_R^{2}}/{{\bf c}_t^{2}}>0$, $\varkappa_{l}^2=1-{{\bf c}_R^{2}}/{{\bf c}_l^{2}}>0$. Notice that $\varkappa_{t}$ and $\varkappa_{l}$ depend on $\sigma$ only.
\medskip

The following inequality is a consequence of the Rayleigh equation (see \cite{LR} and \cite[Ch. 3]{LL}) and will be used in the sequel:
\begin{equation}\label{Req}
16-24B+8B^2-B^3=16\,\frac{{\bf c}_t^2-{\bf c}_R^2}{{\bf c}_l^2}>0,\quad \mbox{where}\quad B:=1-\varkappa_t^{2}.
\end{equation}
\smallskip

\section{Transformation of the problem (\ref{2dim11})}\label{trans}

\begin{figure}[h!]
	\begin{center}
			\includegraphics[width=0.7\linewidth]{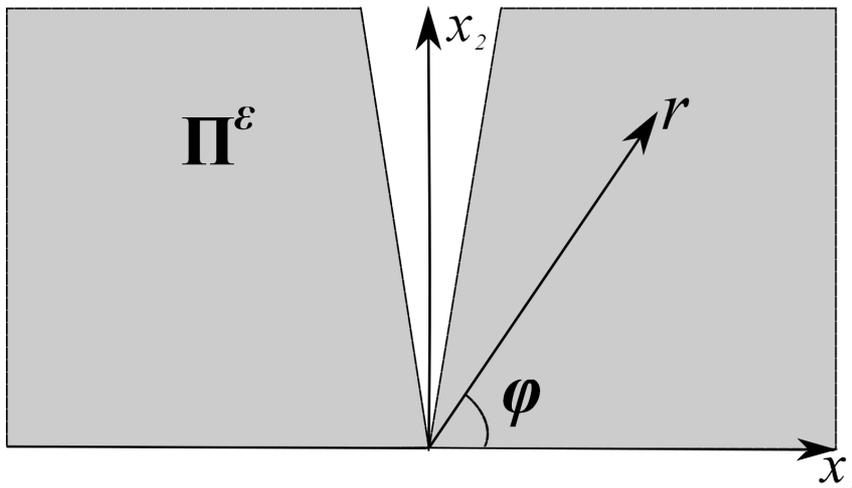}
			\caption{Half-plane with a wedge-like notch 
			} 
			\label{Fig3} 
	\end{center} 
\end{figure}

Using the idea of \cite{SN2011} we cut the angle $\Omega^\e$ along the line $x_1=0$ and transform it by rotations into a half-plane with a wedge-like notch of aperture $2\alpha$ (cf. Fig. \ref{Fig1} and Fig. \ref{Fig3})
\begin{equation*}\label{0}
\Pi^{\e}=\left\{(x_1,x_2)\in \mathbb{R}^2:\ x_2>0,\ |x_1|> \e x_2\right\}.
\end{equation*}

By Remark \ref{rotations}, the equations and the boundary conditions (\ref{2dim11}) are rewritten as follows:
\begin{equation}\label{1}
\aligned
&{\cal L}(\partial_{1},\partial_{2},ik)u^\e  =  \omega^2 u^\e\quad\mbox{in}\quad \Pi^{\e},\\
&{\cal N}(\partial_{1},\partial_{2},ik)u^\e \equiv (\sigma_{12}[u^\e],\sigma_{22}[u^\e],\sigma_{32}[u^\e])(x_1,0)  =   0,
\endaligned
\end{equation}
and should be supplemented with the transmission conditions on the surfaces of the removed infinite sector
\begin{equation}\label{3a}
\aligned
& u^{\e}_{n}(\e x_2,x_2)\nu_{n}^{\e+}+u^{\e}_{n}(-\e x_2,x_2)\nu_{n}^{\e-}=0,\\ 
& u^{\e}_{n}(\e x_2,x_2)\tau_{n}^{\e+}+u^{\e}_{n}(-\e x_2,x_2)\tau_{n}^{\e-}=0,
\\
& u^{\e}_{3}(\e x_2,x_2)=u^{\e}_{3}(-\e x_2,x_2);
\endaligned
\end{equation}
\begin{equation}\label{4a}
\sigma_{nm}[u^\e](\e x_2,x_2)\nu_{m}^{\e+}+ \sigma_{nm}[u^\e](-\e x_2,x_2)\nu_{m}^{\e-}=0,\quad n=1,2,3,
\end{equation}
where $x_2 >0$ and
\begin{equation}\label{nu}
\nu^{\e\pm}=\left(\mp \frac {1}{\sqrt{1+\e^{2}}},\frac {\e}{\sqrt{1+\e^{2}}},0\right) \quad\mbox{and}\quad
\tau^{\e\pm}=\left(-\frac {\e}{\sqrt{1+\e^{2}}},\mp \frac {1}{\sqrt{1+\e^{2}}},0\right)
\end{equation}
are the outward normal and tangential unit vectors, respectively.
Relations (\ref{3a}) and (\ref{4a}) ensure the smoothness of displacements after returning to the original domain $\Omega^\e$.

The variational (weak) formulation of the problem (\ref{1})--(\ref{4a}) refers to the integral identity
\begin{equation}\label{tozhd}
a_{\Pi^\e}(ik;u^{\e},v):=\int\limits_{\Pi^\e}\sigma_{nm}[u^{\e}]\overline{\partial_{m}v_{n}}\,dx_1dx_2=\omega^2 (u^\e,v)_{\Pi^\e},\quad v\in \left(H^{1}_{\sharp}(\Pi^{\e})\right)^3, 
\end{equation}
where bar stands for complex conjugation and $\left(\cdot,\cdot\right)_{\Pi^{\e}}$ for the natural inner product in $\left(L^2 (\Pi^{\e})\right)^3$, while $\left(H^{1}_{\sharp}(\Pi^{\e})\right)^3$ is a subspace of 
vector functions in $\left(H^{1}(\Pi^{\e})\right)^3$ satisfying the condition (\ref{3a}). The quadratic form $a_{\Pi^{\e}}(ik;u,u)$ in $\left(L^2 (\Pi^{\e})\right)^3$ with domain $\big(H^{1}_{\sharp}(\Pi^\e)\big)^3$ is 
symmetric, positive and closed. Therefore, the boundary value problem (\ref{1})--(\ref{4a}) is associated with an unbounded positive definite self-adjoint operator ${\cal A}^\e (ik)$ in the space 
$\left(L^2 (\Pi^{\e})\right)^3$  (see, e.g., \cite[Ch. 10]{BS}). 

It is shown in \cite{K} that the essential spectrum of this operator coincides with the ray $[{\bf c}_R^{2} k^2,+\infty)$. Consequently, only discrete spectrum may occur below the cutoff point 
$\omega_R^{2}={\bf c}_R^{2} k^2$. If this spectrum is non-empty, the corresponding eigenvector generates a localized wedge mode propagating along the edge with a velocity less than ${\bf c}_R$.

\section{The inner expansion}\label{inner}

According to \cite[Ch. 10]{MNP}, formation of a notch (Fig. \ref{Fig3}) in the half-plane 
$$
\Omega^{0}=\mathbb{R}^{2}_+=\{(x_1,x_2)\in\mathbb{R}^2\ :\ x_2>0\}
$$ 
should be looked upon as a regular perturbation due to special transmission conditions 
(\ref{3a}) and (\ref{4a}). Hence, it seems natural to assume that an eigenvalue $\omega^2$ lies near the cutoff point $\omega_R^{2}$ and to choose the following ansatz for the ``inner'' asymptotic expansion of the corresponding eigenvector in the bounded vicinity of the cut:
\begin{equation}\label{ansatz}
u^{\e}=w^{0}+\e w^{1}+\e^{2} w^{2}+\ldots,
\end{equation}
where higher-order terms are not needed in our formal analysis.

As the main asymptotic term $w^{0}$ we take the solution decaying as $x_2\to+\infty$ of the wave propagation problem along the traction-free boundary of an elastic isotropic half-space $\mathbb{R}^{3}_+$, namely 
the Rayleigh wave \cite{LR}, see also \cite[Ch. 3]{LL},
\begin{equation}\label{Rayleigh}
u^R(x_2)=\left(0, u_{2}^{R}(x_2), u_{3}^{R}(x_2)\right),
\end{equation}
\begin{equation}\label{R1}
\aligned
& u^R_{2}(x_2)=iA\varkappa_{l}\Big(e^{-k\varkappa_{l}x_2}-\frac 2{1+\varkappa_{t}^{2}}\,e^{-k\varkappa_{t}x_2}\Big),\\
& u^R_{3}(x_2)=A\Big(e^{-k\varkappa_{l}x_2}-\frac{2\varkappa_{l}\varkappa_{t}}{1+\varkappa_{t}^2}\,e^{-k\varkappa_{t}x_2}\Big).
\endaligned
\end{equation}
Without loss of generality we set $A=1$.

We insert formula (\ref{ansatz}) with $w^0=u^R$ into the transmission conditions (\ref{3a}), (\ref{4a}) and collect terms of order $\e$. As a result, the first correction term $w^1$ should satisfy the problem (here $n=1,2,3$)
\begin{equation}\label{5}
{\cal L}(\partial_{1},\partial_{2},ik)w^{1}=\omega^{2}_{R}w^{1} \quad \mbox{in}\quad\Pi^{0},\qquad
{\cal N}(\partial_{1},\partial_{2},ik)w^{1}(x_1,0)= 0;
\end{equation}
\begin{equation}\label{8}
\aligned
-w^{1}_{n}(+0,x_2)+w^{1}_{n}(-0,x_2) & =
-2u^R_{2}(x_2)\delta_{n,1},
\\
-\sigma_{n1}[w^{1}](+0,x_2)+\sigma_{n1}[w^{1}](-0,x_2) & = -2\sigma_{n2}[u^R](0,x_2),
\endaligned
\qquad x_2>0.
\end{equation}

\subsection{Operator pencil generated by the problem (\ref{5}), (\ref{8})}\label{Ss:pencil}

First, we notice that the homogeneous problem  (\ref{5}), (\ref{8}) (with the right-hand sides in (\ref{8}) equal to zero) appears to be a problem in the half-plane $\mathbb{R}^2_+$. We apply the Fourier--Laplace transform (the Fourier transform without imaginary unit in the exponent)
$$
\widehat w(\xi,x_2)=\frac 1{\sqrt{2\pi}}\int\limits_{\mathbb R}e^{\xi x_1} w(x_1,x_2)\,dx_1
$$
with the complex parameter $\xi$ in a strip
$$
\Sigma_\beta=\{\xi\in\mathbb{C}: |\Re \xi|\le\beta\}
$$
of small width $2\beta>0$ (in what follows we do not display the dependence of $\widehat w$ on $\xi$). This gives a boundary value problem for a system of ordinary differential equations (recall that $\omega_{R}={\bf c}_Rk$)
\begin{equation*}\label{ODE}
%\aligned
{\cal L}(\xi,\partial_2,ik)\widehat w(x_2)-\omega^2_{R} \widehat w(x_2)=0, \quad x_2>0;\qquad
{\cal N}(\xi,\partial_2,ik)\widehat w(0)=0.
\end{equation*}
Here ${\cal L}(\xi,\partial_{2},ik)$ is the matrix differential operator
\begin{equation*}
\begin{bmatrix} \mu(k^2-\partial_{2}^{2})-(\lambda+2\mu)\xi^{2} & -(\lambda+\mu)\xi\partial_{2} & -(\lambda+\mu)i\xi k \\ -(\lambda+\mu)\xi\partial_{2} & -(\lambda+2\mu)\partial^{2}_{2}+\mu(k^2-\xi^{2})& -
(\lambda+\mu)i k\partial_{2} \\ -(\lambda+\mu)i\xi k & -(\lambda+\mu)ik\partial_{2} & -\mu(\xi^{2}+\partial^{2}_{2} )+(\lambda+2\mu)k^2\end{bmatrix},
\end{equation*}
and analogously
\begin{equation*}
{\cal N}(\xi,\partial_{2},ik)=
\begin{bmatrix} -\mu\partial_2 & -\mu \xi & 0 \\ -\lambda \xi & -(\lambda+2\mu)\partial_2 & -\lambda ik \\ 0 & -\mu ik & -\mu \partial_2 \end{bmatrix}.
\end{equation*}

\begin{Theorem}\label{L1}
There is a $\beta>0$ depending only on $\lambda$, $\mu$ and $k$ such that the spectrum of the (quadratic) operator pencil 
\begin{equation}\label{pencil}
\aligned
\Sigma_\beta \ni \xi \mapsto \mathfrak{A}(\xi) &\,:\ 
\big(H^2(\mathbb{R}_+)\big)^3 \to \big(L^2(\mathbb{R}_+)\big)^3\times \mathbb{C}^3
\\
\mathfrak{A}(\xi) &\, := \left({\cal L}(\xi,\partial_2,ik)-\omega_R^2; \, {\cal N}(\xi,\partial_2,ik)\big|_{x_2=0}\right)
\endaligned
\end{equation}
in the strip $\Sigma_\beta$ consists of a unique eigenvalue $\xi=0$ of algebraic multiplicity two. Corresponding Jordan block consists of eigenvector $\widehat w^0(x_2)=u^R(x_2)$ and associated vector 
\begin{equation}\label{adjoint}
\widehat w^1(x_2)=\big(-ik^{-1}u^R_{3}(x_2),0,0\big).
\end{equation}
\end{Theorem}

\noindent{\bf Proof.} We proceed in four steps.

{\bf Step 1}. We claim that imaginary axis apart from the origin is free of the spectrum of the pencil (\ref{pencil}). Indeed,
let $\xi=i\eta$, $\eta\in \mathbb{R}$. 
Then 
\begin{equation}\label{identity}
\aligned
{\cal L}(i\eta,\partial_2,ik)\equiv M^{-1}{\cal L}(0,\partial_2,i\sqrt{k^2+\eta^2})M,\\
{\cal N}(i\eta,\partial_2,ik)\equiv M^{-1}{\cal N}(0,\partial_2,i\sqrt{k^2+\eta^2})M,
\endaligned
\end{equation}
where $M$ is the matrix of rotation in the plane $x_1Ox_3$ by the angle $\arctan(\eta/k)$.

By \cite[Lemma 2(b)]{BDJ}, the identity (\ref{identity}) gives
\begin{equation}\label{positive1}
\int\limits_0^\infty {\cal L}(i\eta,\partial_2,ik) {\mathfrak u}(x_2)\cdot {\mathfrak u}(x_2)\,dx_2 \ge {\bf c}_R^2(k^2+\eta^2)\|{\mathfrak u};\big(L^2(\mathbb{R}_+)\big)^3\|^2
\end{equation}
for arbitrary vector function ${\mathfrak u}\in \big(H^2(\mathbb{R}_+)\big)^3$ which satisfies the boundary condition ${\cal N}(i\eta,\partial_2,ik){\mathfrak u}(0)=0$. Here and further the central dot stands for the scalar product in $\mathbb C^3$. 

In the language of wave propagation, the inequality (\ref{positive1}) means that any eigenvector of the operator
\begin{equation*}%\label{pencil1}
\left({\cal L}(i\eta,\partial_2,ik)-\omega^2;\  {\cal N}(i\eta,\partial_2,ik)\right)
\end{equation*}
generates a surface wave in the half-space running with the velocity $\omega/\sqrt{k^2+\eta^2}$ which cannot be less than ${\bf c}_R$.

In particular, (\ref{positive1}) implies that an eigenvector $\widehat w^0(x_2)$ of the pencil (\ref{pencil}) with $\xi=i\eta$, $\eta\in\mathbb{R}$, can exist only for $\eta=0$, and the claim follows.
\smallskip

{\bf Step 2}. Now we study the zero eigenvalue of the pencil (\ref{pencil}). For the corresponding eigenvector we have
\begin{equation}\label{ODE1}
{\cal L}(0,\partial_2,ik)\widehat w^0(x_2)=\omega^2_{R} \widehat w^0(x_2),\quad x_2>0;\qquad
{\cal N}(0,\partial_2,ik)\widehat w^0(0)=0,
\end{equation}
and thus $\widehat w^0(x_2)=u^R(x_2)$ (up to a multiplicative constant).

For the associated vector of rank 1, we obtain the problem
\begin{equation}\label{ODE2}
\aligned
{\cal L}(0,\partial_2,ik)\widehat w^1(x_2)-\omega^2_{R} \widehat w^1(x_2)=-\partial_\xi {\cal L} (0,\partial_{2},ik)\widehat w^0(x_2),&& x_2>0;\\
{\cal N}(0,\partial_2,ik)\widehat w^1(0)=-\partial_\xi {\cal N}(0,\partial_{2},ik)\widehat w^0(0).
\endaligned
\end{equation}
By direct calculation we check that the vector function (\ref{adjoint}) solves this problem.

Finally, the associated vector of rank 2 should be a solution of the problem
\begin{equation}\label{ODE3}
\aligned
{\cal L}(0,\partial_2,ik)\widehat w^2(x_2)\,-&\, \omega^2_{R} \widehat w^2(x_2)\\
=&-\partial_\xi {\cal L}(0,\partial_{2},ik)\widehat w^1(x_2)-\tfrac 12\partial^2_\xi {\cal L}(0,\partial_{2},ik)\widehat w^0(x_2),\quad x_2>0;\\
{\cal N}(0,\partial_2,ik)\widehat w^2(0)= &-\partial_\xi {\cal N}(0,\partial_{2},ik)\widehat w^1(0)-\tfrac 12\partial^2_\xi {\cal N}(0,\partial_{2},ik)\widehat w^0(0).
\endaligned
\end{equation}
Assume that such solution exists. We multiply this equation by the eigenvector $\widehat w^0$ and integrate over $\mathbb{R}_+$. Then integrating by parts and using (\ref{ODE1}) yield
\begin{equation}\label{b=0}
\aligned
 0=b&:= -\,\frac{\lambda}{ik}\,u^R_{3}(0)\overline{u_{2}^{R}}(0)\\
&+\int\limits_0^\infty
\big(i(\lambda+\mu)u^{R\prime}_3(x_2)\overline{u_{2}^{R}}(x_2)-(\lambda+2\mu)|u_{3}^{R}(x_2)|^{2}-\mu|u^R_2(x_2)|^2\big)dx_2.
\endaligned
 \end{equation}
A direct calculation based on the explicit formulae (\ref{R1}) gives
\begin{equation}
b =-(\lambda+2\mu)\,\frac{B^{3}(16-24B+8B^2-B^3)((1-B)^{2}(7-2B)+1)}{128k(1-B)^{\frac{5}{2}}(2-B)^{2}},
\label{b}
 \end{equation}
we recall that $B=1-\varkappa_t^{2}\in(0,1)$. Since $b<0$ according to (\ref{Req}), the relation (\ref{b=0}) is impossible. Thus, the algebraic multiplicity of the zero eigenvalue is just $2$.
\smallskip

{\bf Step 3}. By \cite[Lemma 2(a)]{BDJ}, the inequality
$$
\int\limits_0^\infty {\cal L}(0,\partial_2,ik) {\mathfrak u}(x_2)\cdot {\mathfrak u}(x_2)\,dx_2 \ge {\bf c}_t^2k^2\|{\mathfrak u};\big(L^2(\mathbb{R}_+)\big)^3\|^2
$$
holds for all vector functions ${\mathfrak u}\in \big(H^2(\mathbb{R}_+)\big)^3$ which satisfy the boundary condition ${\cal N}(0,\partial_2,ik){\mathfrak u}(0)=0$ and are
orthogonal to $\widehat w^0$ in $\big(L^2(\mathbb{R}_+)\big)^3$. Since $\omega_R={\bf c}_Rk<{\bf c}_tk$, this implies that $0$ is an isolated simple eigenvalue of the (self-adjoint) operator $\mathfrak{A}(0)$. 

To study the operators $\mathfrak{A}(\xi)$ for $\xi$ close to zero we notice that $\mathfrak{A}(0)$ is positive semidefinite by (\ref{positive1}). Thus, the operator
$$
\widetilde{\mathfrak{A}}(\xi) := \left({\cal L}(\xi,\partial_2,ik)-\omega_R^2+1; \, {\cal N}(\xi,\partial_2,ik)\big|_{x_2=0}\right)
$$
is invertible if $|\xi|$ is sufficiently small, and the operator
$$
\widehat{\mathfrak{A}}(\xi)=\big(\widetilde{\mathfrak{A}}(\xi)\big)^{-1}{\mathfrak{A}}(\xi)
$$
is bounded in $\big(H^2(\mathbb{R}_+)\big)^3$. Moreover, $0$ is an isolated simple eigenvalue of the operator $\widehat{\mathfrak{A}}(0)$. 

Theorem on stability of the root multiplicities, see \cite[Ch. I, Theorem 3.1]{GK}, shows that there is a $\delta>0$ such that if $|\xi|<\delta$ then the circle $|z|<\delta$ does not contain points of essential spectrum of the operator $\widehat{\mathfrak{A}}(\xi)$ and contains a unique simple eigenvalue $\widehat z(\xi)$. The same holds for the original operator ${\mathfrak{A}}(\xi)$ since  
$\widetilde{\mathfrak{A}}(\xi)$ is an isomorphism.

The results of {\bf Step 2} imply that 
$$
\frac {d\widehat z}{d\xi}(0)=0 \qquad\mbox{and}\qquad \frac {d^2\widehat z}{d\xi^2}(0)\ne0,
$$
cf. (\ref{CC}). Therefore, $\widehat z(\xi)\ne0$ if $|\xi|\ne0$ is sufficiently small, and thus a neighborhood of the origin contains a unique eigenvalue $\xi=0$ of the pencil (\ref{pencil}).
\smallskip

{\bf Step 4}. By (\ref{positive1}), the self-adjoint operator $\mathfrak{A}(i\eta)$  with $\eta\in\mathbb{R}\setminus\{0\}$ is invertible. Thus, it possesses the isomorphism property which is preserved under small perturbations. Therefore, a neighborhood of imaginary axis contains a unique eigenvalue $\xi=0$ of the pencil (\ref{pencil}). It remains to show that this neighborhood cannot narrow at infinity. 

Integration by parts gives
\begin{equation}\label{identity1}
\aligned
& \int\limits_0^\infty {\cal L}(\xi,\partial_2,ik) {\mathfrak u}(x_2)\cdot {\mathfrak u}(x_2)\,dx_2=
\int\limits_0^\infty \Big[ \mu|{\mathfrak u}_1'|^2+(\mu k^2-(\lambda+2\mu)\xi^2)|{\mathfrak u}_1|^2\\
+ &\ (\lambda+2\mu)|{\mathfrak u}_2'|^2 +\mu(k^2-\xi^2)|{\mathfrak u}_2|^2+\mu|{\mathfrak u}_3'|^2 +((\lambda+2\mu) k^2-\mu\xi^2)|{\mathfrak u}_3|^2\\
+ &\ \lambda\xi({\mathfrak u}_1\overline{\mathfrak u}_2' -{\mathfrak u}_2'\overline{\mathfrak u}_1)
+\mu\xi({\mathfrak u}_2\overline{\mathfrak u}_1' -{\mathfrak u}_1'\overline{\mathfrak u}_2) +\lambda ik({\mathfrak u}_3\overline{\mathfrak u}_2' -{\mathfrak u}_2'\overline{\mathfrak u}_3) 
+\mu ik({\mathfrak u}_2\overline{\mathfrak u}_3' -{\mathfrak u}_3'\overline{\mathfrak u}_2)\\
- &\ (\lambda+\mu)ik\xi({\mathfrak u}_1\overline{\mathfrak u}_3 +{\mathfrak u}_3\overline{\mathfrak u}_1)\Big]dx_2=: {\mathfrak a}(\xi,ik;{\mathfrak u},{\mathfrak u})
\endaligned
\end{equation}
for arbitrary vector function ${\mathfrak u}\in \big(H^2(\mathbb{R}_+)\big)^3$ which satisfies the boundary condition ${\cal N}(\xi,\partial_2,ik){\mathfrak u}(0)=0$.

Let $\xi=\zeta+i\eta$. Then for arbitrary vector function ${\mathfrak u}\in \big(H^1(\mathbb{R}_+)\big)^3$ we estimate
\begin{equation}\label{estimate}
\aligned
& |{\mathfrak a}(\xi,ik;{\mathfrak u},{\mathfrak u}) -{\mathfrak a}(i\eta,ik;{\mathfrak u},{\mathfrak u})|\le
C|\zeta|\,\|{\mathfrak u};\big(L^2(\mathbb{R}_+)\big)^3\| \\ 
\times & \big(\|{\mathfrak u}';\big(L^2(\mathbb{R}_+)\big)^3\|+
(k+|\zeta|+|\eta|)\|{\mathfrak u};\big(L^2(\mathbb{R}_+)\big)^3\|\big),
\endaligned
\end{equation}
where $C$ depends on $\lambda$ and $\mu$ only. Furthermore, 
\begin{equation*}
\aligned
{\mathfrak a}(i\eta,ik;{\mathfrak u},{\mathfrak u})=&\int\limits_0^\infty \Big[ \mu\big(|{\mathfrak v}_1'|^2+(k^2+\eta^2)|{\mathfrak v}_1|^2\big) \\
+&\ (\lambda+2\mu)\big(|{\mathfrak v}_2'|^2+(k^2+\eta^2)|{\mathfrak v}_3|^2\big)+2\lambda\sqrt{k^2+\eta^2}\Im\big( {\mathfrak v}_2' \overline{\mathfrak v}_3\big) \\
+&\ \mu\big(|{\mathfrak v}_3'|^2+(k^2+\eta^2)|{\mathfrak v}_2|^2+2\sqrt{k^2+\eta^2}\Im\big( {\mathfrak v}_3' \overline{\mathfrak v}_2\big)\big)\Big] dx_2.
\endaligned
\end{equation*}
Here, ${\mathfrak v}=M{\mathfrak u}$ and the matrix $M$ was introduced in (\ref{identity}).
Using the inequality $\lambda+\frac{2}{3}\mu>0$, we estimate
\begin{equation}\label{positive2}
\aligned
{\mathfrak a}(i\eta,ik;{\mathfrak u},{\mathfrak u})\ge&\ \mu \int\limits_0^\infty \Big[ \big(|{\mathfrak v}_1'|^2+(k^2+\eta^2)|{\mathfrak v}_1|^2\big) \\
+&\ \frac{2}{3}\big(|{\mathfrak v}_2'|^2+(k^2+\eta^2)|{\mathfrak v}_3|^2\big) 
+\gamma\big(|{\mathfrak v}_3'|^2-2(k^2+\eta^2)|{\mathfrak v}_2|^2\big)\Big] dx_2\\
\ge &\ \gamma\mu \Big(\|{\mathfrak u}';\big(L^2(\mathbb{R}_+)\big)^3\|^2-2(k^2+\eta^2) \|{\mathfrak u};\big(L^2(\mathbb{R}_+)\big)^3\|^2\Big)
\endaligned
\end{equation}
with arbitrary $\gamma<\frac 12$.

From (\ref{positive1}) and (\ref{identity1}) we derive
\begin{equation}\label{positive3}
{\mathfrak a}(i\eta,ik;{\mathfrak u},{\mathfrak u})\ge {\bf c}_R^2(k^2+\eta^2)\|{\mathfrak u};\big(L^2(\mathbb{R}_+)\big)^3\|^2
\end{equation}
for arbitrary vector function ${\mathfrak u}\in \big(H^2(\mathbb{R}_+)\big)^3$ which satisfies the boundary condition ${\cal N}(\xi,\partial_2,ik){\mathfrak u}(0)=0$. By approximation, this inequality holds for all ${\mathfrak u}\in \big(H^1(\mathbb{R}_+)\big)^3$. We sum up (\ref{positive2}) and (\ref{positive3}), choose $\gamma={\bf c}_R^2/(3\mu)$ and conclude that
\begin{equation}\label{positive4}
{\mathfrak a}(i\eta,ik;{\mathfrak u},{\mathfrak u})\ge \frac {{\bf c}_R^2}6\Big(\|{\mathfrak u}';\big(L^2(\mathbb{R}_+)\big)^3\|^2+(k^2+\eta^2)\|{\mathfrak u};\big(L^2(\mathbb{R}_+)\big)^3\|^2\Big).
\end{equation}
The relations (\ref{estimate}) and (\ref{positive4}) ensure that for any ${\mathfrak u}\in \big(H^1(\mathbb{R}_+)\big)^3$ we have
\begin{equation}\label{positive5}
\aligned
& \left|{\mathfrak a}(\xi,ik;{\mathfrak u},{\mathfrak u}) -\omega_R^2\|{\mathfrak u};\big(L^2(\mathbb{R}_+)\big)^3\|^2\right| \ge \Big(\frac {{\bf c}_R^2}6-C\epsilon\Big)\|
{\mathfrak u}';\big(L^2(\mathbb{R}_+)\big)^3\|^2\\
 + &\ \Big(\frac {{\bf c}_R^2}6 \eta^2-\frac {5{\bf c}_R^2}6 k^2 -\frac {C\zeta^2}{4\epsilon}-C|\zeta|(k+|\zeta|+|\eta|)\Big)\|{\mathfrak u};\big(L^2(\mathbb{R}_+)\big)^3\|^2.
\endaligned
\end{equation}
The right-hand side of (\ref{positive5}) is bounded away from zero for any given $k$ and $\zeta$ if $\epsilon>0$ is small and $|\eta|$ is large. By (\ref{identity1}), the operator $\mathfrak{A}(\xi)$  with such $\xi=\zeta+i\eta$ and its formally adjoint operator  $\mathfrak{A}(-\overline{\xi})$ are invertible. Thus, they both possess the isomorphism property. Therefore, a neighborhood of imaginary axis considered at the beginning of {\bf Step 4} contains a strip. The proof is completed.\hfill$\square$

\begin{Remark}\label{vishik}
Notice that by the inequality (\ref{positive3}) a solution of the problem 
$$
\mathfrak{A}(i\eta){\mathfrak u}=({\bf f,g})\in \big(L^2(\mathbb{R}_+)\big)^3\times \mathbb{C}^3
$$
with $\eta\in\mathbb{R}$, $|\eta|\ge\eta_0>0$ satisfies the following estimate of a parameter-dependent norm uniform in $\eta$: 
$$
\|{\mathfrak u}'';\big(L^2(\mathbb{R}_+)\big)^3\|^2
+\eta^2\|{\mathfrak u};\big(L^2(\mathbb{R}_+)\big)^3\|^2\le c(\eta_0)\Big(\|{\bf f};\big(L^2(\mathbb{R}_+)\big)^3\|^2+|{\bf g}|^2\Big).
$$
\end{Remark}

In particular, Theorem \ref{L1} implies that a general solution of the homogeneous problem (\ref{5}), (\ref{8}) with polynomial growth in $x_1$ is a linear combination $C_u u^R(x_2)+C_v v^R(x_1,x_2)$, 
where $u^R$ is the Rayleigh wave (\ref{Rayleigh}) while 
\begin{equation*}
v^{R}(x)=x_1\widehat w^0(x_2)+\widehat w^1(x_2)=\big(-ik^{-1}u^R_{3}(x_2), x_1u_{2}^{R}(x_2), x_1u_{3}^{R}(x_2)\big),
\end{equation*}
see also \cite{APK}.

\subsection{Application of the Kondrat'ev theory}\label{Kondr}

Let us comment on certain peculiarities in application of the Kondrat'ev theory \cite{Ko} (see also the pioneering paper \cite{AgNi}, the monograph \cite[Ch. 3 and 5]{NaPl} and the survey paper \cite[\S3]{na262}) 
to assure that any solution of the non-homogeneous problem (\ref{5}), (\ref{8}) with the polynomial growth as $x_1\to\pm\infty$ gets the form 
\begin{equation}\label{S1}
w^1(x)=\sum_{\pm}\chi(\pm x_1)(c_u^{\pm}u^R(x_2)+c_v^{\pm}v^R(x))+{\widetilde{w}^{1}}(x).
\end{equation}
 Here, $\chi$ is a smooth cutoff function 
\begin{equation}\label{chi}
 0\le\chi\le1,\qquad\chi(t)=1 \quad\mbox{for}\quad t>2, \qquad \chi(t)=0 \quad\mbox{for}\quad t<1, 
\end{equation}
 and the remainder $\widetilde{w}^{1}(x)$ decays exponentially as $|x|\to\infty$.

\subsubsection{An auxiliary problem in the half-plane}

First of all, we study the problem in the half-plane
\begin{equation}\label{he1}
{\cal L}(\partial_{1},\partial_{2},ik)\mathfrak{w}-\omega^{2}_{R}\mathfrak{w}=\mathfrak{f}\quad\mbox{in}\quad\mathbb{R}^2_+,\qquad{\cal N}(\partial_{1},\partial_{2},ik)\mathfrak{w}=\mathfrak{g}
\quad\mbox{on}\quad\partial\mathbb{R}^2_+. 
\end{equation}

To deploy the corresponding technique, we ought to regard the half-plane as the cylinder  $\mathbb{R}\times\mathbb{R}_+$. A non-standard detail is that the cross-section $\mathbb{R}_+\ni x_2$ is an unbounded domain 
in contrast to the traditional formulation of elliptic boundary value problems in cylinders with bounded cross-sections. However, this distortion is compensated by Theorem \ref{L1} which provides the pencil 
(\ref{pencil}) with all the required properties\footnote{Cf. \cite[Proposition 3.2(2)]{na262} for the Neumann problem for a formally self-adjoint elliptic system in cylinders with bounded cross-sections.}. 
Below we outline the customary steps to conclude on the existence of a solution to (\ref{he1}) and its decomposition (\ref{S1}) (the detailed description can be found, e.g., in \cite[Ch. 3, \S\S1-2]{NaPl}).

We introduce the weighted Lebesgue (for $\ell=0$) and Sobolev (for $\ell\in\mathbb{N}$) spaces $W^{\ell}_{\beta}(\mathbb{R}^2_+)$ with the norms
$$
 \|\mathfrak{w};W^{\ell}_{\beta}(\mathbb{R}^2_+)\|=\|e^{\beta x_1}\mathfrak{w};H^{\ell}(\mathbb{R}^2_+)\|.
$$
Here the weight index $\beta>0$ is defined in Theorem \ref{L1}. Note that a function $\mathfrak{w}\in W^{\ell}_{\beta}(\mathbb{R}^2_+)$ decays exponentially as $x_1\to+\infty$ but may grow as $x_1\to-\infty$. 
By $W^{\ell+\frac 12}_{\beta}(\partial\mathbb{R}^2_+)$ we denote the weighted Sobolev--Slobodetskii space supplied with the natural trace norm
$$
 \|\mathfrak{g};W^{\ell+\frac 12}_{\beta}(\partial\mathbb{R}^2_+)\|=\inf\big\{\|G;W^{\ell}_{\beta}(\mathbb{R}^2_+)\|\ :\ G\big|_{\partial\mathbb{R}^2_+}=\mathfrak{g}\big\}.
$$

We consider the problem (\ref{he1}) with the right-hand sides
\begin{equation}\label{he2}
 \mathfrak{f}\in \big(W^{\ell}_{\beta}(\mathbb{R}^2_+)\big)^3\cap
\big(W^{\ell}_{-\beta}(\mathbb{R}^2_+)\big)^3,\qquad
\mathfrak{g}\in \big(W^{\ell+\frac 12}_{\beta}(\partial\mathbb{R}^2_+)\big)^3\cap \big(W^{\ell+\frac 12}_{-\beta}(\partial\mathbb{R}^2_+)\big)^3,
\end{equation}
which decay as $x_1\to\pm\infty$. In this way, the problem (\ref{he1}) is associated with two operators
\begin{equation}\label{he3}
A_{\pm}:\ \big(W^{\ell+2}_{\pm\beta}(\mathbb{R}^2_+)\big)^3\to \big(W^{\ell}_{\pm\beta}(\mathbb{R}^2_+)\big)^3\times \big(W^{\ell+\frac 12}_{\pm\beta}(\partial\mathbb{R}^2_+)\big)^3.
\end{equation}

By Theorem \ref{L1}, the lines $I_{\pm}=\{\xi\in\mathbb{C}:\Re\xi=\pm\beta\}$ are free of the spectrum of the pencil (\ref{pencil}). Therefore, both the operators (\ref{he3}) are isomorphisms. This conclusion is made 
by the standard scheme \cite{AgNi}, \cite{Ko}: the direct Fourier transform, the resolvent ${\cal R}(\xi)=(\mathfrak{A}(\xi))^{-1}$ with $\xi\in I_{\pm}$ and the inverse Fourier transform\footnote{In the cited literature the 
pencil always possesses the important property: for any $\xi_1,\xi_2\in\mathbb{C}$ the difference $\mathfrak{A}(\xi_1)-\mathfrak{A}(\xi_2)$ is a compact operator. In our case this property is absent and this was the 
very reason to prove Theorem \ref{L1} individually.}. Therefore, in view of (\ref{he2}) we obtain two solutions $\mathfrak{w}^{\pm}\in \big(W^{\ell+2}_{\pm\beta}(\mathbb{R}^2_+)\big)^3$ of the problem (\ref{he1}). 
These solutions are given by integrals along the lines $I_{\pm}$. By Theorem \ref{L1}, $\xi\mapsto{\cal R}(\xi)$ is a meromorphic operator function in the strip $\Sigma_\beta$ between the lines $I_-$ and $I_+$ with the only pole 
of rank 2 at the point $\xi=0$. Therefore, the Cauchy residual theorem gives the relation
\begin{equation}\label{he4}
\mathfrak{w}^{+}(x_1,x_2)-\mathfrak{w}^{-}(x_1,x_2)=a_{u}u^{R}(x_2)+a_{v}v^{R}(x),
\end{equation}
and the following estimate is valid:
$$
\sum\limits_{\pm}\big\|\mathfrak{w}^{\pm};\big(W^{\ell+2}_{\pm\beta}(\mathbb{R}^2_+)\big)^3\big\|+|a_{u}|+|a_{v}|\leq C\sum\limits_{\pm}\Big(\big\|\mathfrak{f};\big(W^{\ell}_{\pm\beta}(\mathbb{R}^2_+)\big)^3\big\|
+\big(\big\|\mathfrak{g};W^{\ell+\frac 12}_{\pm\beta}(\partial\mathbb{R}^2_+)\big)^3\big\|\Big),
$$
where $C$ depends only on $\lambda$, $\mu$ and $k$.

Since $\mathfrak{w}^+$ decays exponentially as $x_1\to+\infty$ and $\mathfrak{w}^-$ does as $x_1\to-\infty$, formula (\ref{he4}) implies that they both take the form (\ref{S1}) (with $c_u^+=c_v^+=0$ and $c_u^-=c_v^-=0$, 
respectively). At the same time, for arbitrary constants $C_u$ and $C_v$ the linear combination $\mathfrak{w}^++C_u u^R+C_v v^R$ still verifies the problem (\ref{he1}) and takes the form (\ref{S1}).

\begin{Remark}\label{remain}
The remainder term $\widetilde{w}^1$ in (\ref{S1}) decays exponentially as $|x_1|\to\infty$ but in general not as $x_2\to+\infty$. However, it is easy to see that if $\mathfrak{f}$ decays exponentially as $x_2\to+\infty$ 
then $\widetilde{w}^1$ does the same. 
\end{Remark}

\subsubsection{Reduction of the problem (\ref{5}), (\ref{8}) to the form (\ref{he1})}

We proceed in two steps. According to (\ref{R1}), the right-hand sides in (\ref{8}) are smooth and decay exponentially as $x_2\to+\infty$. Therefore, there is a smooth vector function $w^{\sharp}$ in the first quadrant 
which also decays exponentially as $x_2\to+\infty$ and meets the boundary conditions
$$
w^{\sharp}_{n}(+0,x_2) 
=u^R_{2}(x_2)\delta_{n,1},\quad 
\sigma_{n1}[w^{\sharp}](+0,x_2)
= \sigma_{n2}[u^R](0,x_2),
\qquad x_2>0,\quad n=1,2,3.
$$
We define $w^{\sharp}$ in the second quadrant as
$$
w^{\sharp}_{1}(-x_1, x_2):=-w^{\sharp}_{1}(x_1,x_2), \quad 
w^{\sharp}_{2}(-x_1, x_2):=w^{\sharp}_{2}(x_1,x_2), \quad 
w^{\sharp}_{3}(-x_1, x_2):=w^{\sharp}_{3}(x_1,x_2).
$$
Then the composite function
$$
w^{\flat}(x)=
\begin{cases}
w^1(x)-(1-\chi(x_1))w^{\sharp}(x), & x_1>0,\\
w^1(x)-(1-\chi(-x_1))w^{\sharp}(x), & x_1<0
\end{cases}
$$
(recall that the cutoff function $\chi$ is introduced in (\ref{chi})) meets the homogeneous jump conditions (\ref{8}) and therefore satisfies the problem in the half-plane
$$
{\cal L}(\partial_{1},\partial_{2},ik)w^{\flat}-\omega^{2}_{R}w^{\flat}=f^{\flat}\quad\mbox{in}\quad\mathbb{R}^2_+,\qquad{\cal N}(\partial_{1},\partial_{2},ik)w^{\flat}=g^{\flat}\quad\mbox{on}\quad\partial\mathbb{R}^2_+. 
$$
Moreover, the function $f^{\flat}$ is supported in the semi-infinite strip and decays exponentially as $x_2\to+\infty$, while $g^{\flat}$ is supported in the segment $[0,2]$. Thus, 
$f^{\flat}\in\big(W^0_{\pm\beta}(\mathbb{R}^2_+)\big)^3$ and $g^{\flat}\in\big(W^0_{\pm\beta}(\partial\mathbb{R}^2_+)\big)^3$. 

Unfortunately, $g^{\flat}$ has a jump at the origin, $\mathbb{R}^3\ni g^{\bullet} :=g^{\flat}(+0)-g^{\flat}(-0)\ne0$. Therefore, $g^{\flat}\notin \big(W^{\frac 12}_{\pm\beta}(\partial\mathbb{R}^2_+)\big)^3$, and the 
previous argument at $\ell=0$ is not applicable. To come over this obstacle we need the second step.
\medskip

To compensate for the above-mentioned discontinuity, we again employ the Kondrat'ev theory, now in a domain with a corner point at the boundary. Recall that $(r,\phi)$ are the polar coordinates of the point $(x_1,x_2)$, define
$$
\mathfrak{L}(\partial_r, \partial_{\phi}):= {\cal L}(\partial_{1},\partial_{2},0), \qquad \mathfrak{N}(\partial_r, \partial_{\phi}):= {\cal N}(\partial_{1},\partial_{2},0)
$$
and consider the problem
\begin{equation}\label{ugol}
\aligned
& \mathfrak{L}(\partial_r, \partial_{\phi})\Psi=0, \quad (r,\phi)\in \mathbb{R}_+\times(0,\pi);\\
& \mathfrak{N}(\partial_r, \partial_{\phi})\Psi(r,0)=\frac 12 g^{\bullet},\quad 
\mathfrak{N}(\partial_r, \partial_{\phi})\Psi(r,\pi)=-\frac 12 g^{\bullet}.
\endaligned
\end{equation}
Using the Kondrat'ev technicality (see, e.g., \cite[Lemma 3.5.11]{NaPl}) we find out a particular solution of (\ref{ugol}) in the form \begin{equation}\label{he5}
\Psi(r,\phi)=r(\psi^0(\phi)\log(r)+\psi^1(\phi)),
\end{equation}
where $\psi^0,\psi^1\in\big({\cal C}^\infty([0,\pi])\big)^3$.

\begin{Remark}\label{sing-origin}
According to the classical works \cite{Wi} and \cite{Ko}, the displacement field $u^\e$ solving the problem (\ref{2dim11}) in the angle $\Omega^\e$ gets the singular component 
$r^{t(\e)}\Phi^\e((1-\frac {2\alpha}{\pi})\phi+\alpha)$ near the origin. Here $t(\e)$ and $\Phi^\e\in \big({\cal C}^\infty([0,\pi])\big)^3$ depend smoothly on $\e\ge0$, and the following relations hold as $\e\to+0$:
$$
t(\e)=1+\e t^1+O(\e^2), \qquad \Phi^\e((1-\tfrac {2\alpha}{\pi})\phi+\alpha)=\hat\Phi^0(\phi)+\e\hat\Phi^1(\phi)+O(\e^2).
$$
Moreover, $r\hat\Phi^0(\phi)$ is a vector function linear in $x$. We then have
$$
r^{t(\e)}\Phi^\e((1-\tfrac {2\alpha}{\pi})\phi+\alpha)=r\hat\Phi^0(\phi)+\e r
(t^1\hat\Phi^0(\phi)\log(r)+\hat\Phi^1(\phi))+O(\e^2)
$$
that explains the appearance of the logarithmic-dependent term (\ref{he5}) in the asymptotic ansatz.
\end{Remark}

Finally, we set
$$
\mathfrak{w}(x)=w^{\flat}(x)-(1-\chi(r))\Psi(r,\phi)
$$
and claim that $\mathfrak{w}$ meets the problem (\ref{he1}) with right-hand sides satisfying (\ref{he2}) for $\ell=0$.

Indeed, by (\ref{ugol}), the function $\mathfrak{f}-f^{\flat}$ is compactly supported and contains only the first derivatives of $\Psi$. Furthermore,  
$$
\mathfrak{g}(x_1)=g^{\flat}(x_1)-(1-\chi(x_1))\,\frac{\text{sign}(x_1)}2\,g^{\bullet} + h(x_1)\Psi|_{x_2=0},
$$ 
where $h$ is a compactly supported smooth function. Thus the claim follows.

\begin{Remark}\label{Psi}
Notice that for any $\e>0$ we have $u^\e\in \big(H^2_{\rm{loc}}(\Omega^{\e})\big)^3$, see \cite{Wi}, \cite{KM}. In contrast, the vector function $w^1$ does not belong to $\big(H^2_{\rm{loc}}(\mathbb{R}^2_+)\big)^3$ because of singularity of $\Psi$ at the origin. However, $\Psi\in \big(W^{2,p}_{\rm{loc}}(\mathbb{R}^2_+)\big)^3$ for any $p<2$ and therefore $w^1\in \big(W^{2,p}_{\rm{loc}}(\mathbb{R}^2_+)\big)^3$ for any $p<2$ that is sufficient for our purposes.
\end{Remark}

\subsection{The completion of the construction }

Summing up, we have established the solvability of the problem (\ref{5}), (\ref{8}) and a general form (\ref{S1}) of the solution, while the exponential decay of the remainder follows from Remark \ref{remain}. 
To complete the construction of the first correction term in (\ref{ansatz}) we need to compute the constants $c_u^{\pm}$ and $c_v^{\pm}$. Since $w^1$ is defined up to a summand $C_u u^R+C_v v^R$, we can assume 
without loss of generality that $c^+_u + c_u^-=0$ and $c^+_v + c_v^-=0$.

We apply the method \cite{MaPl1}, see also \cite[Ch. 3, \S3]{NaPl}. Namely, we insert $w^1$ and $u^R$ into Green's formula in the rectangles $Q^{\pm}(T)=\{x\in{\mathbb R}^2\ :\ 0<\pm x_1<T,\ 0<x_2<T \}$: 
\begin{multline*}
\int\limits_{\partial Q^{\pm}_T} \big({\cal N}(\partial_{1},\partial_{2},ik)w^{1}\cdot u^R - w^{1}\cdot  {\cal N}(\partial_{1},\partial_{2},ik)u^R\big)\\
=
\int\limits_{Q^{\pm}_T}\big({\cal L}(\partial_{1},\partial_{2},ik)w^{1}\cdot u^R - w^{1}\cdot {\cal L}(\partial_{1},\partial_{2},ik)u^R\big)=0
\end{multline*}
and then send $T$ to $+\infty$. 

Integrals over the line $x_2=0$ equal zero by the traction-free boundary condition while integrals over the line $x_2=T$ vanish by decay of $w^1$ and $u^R$. This gives
\begin{multline*}
\lim\limits_{T\to\infty}\int\limits_{0}^{\infty} \big({\cal N}(\partial_{1},\partial_{2},ik)w^{1}\cdot u^R - w^{1}\cdot  {\cal N}(\partial_{1},\partial_{2},ik)u^R\big)\,dx_2\bigg|_{x_1=-T}^{x_1=T}\\
=\int\limits_{0}^{\infty} \big({\cal N}(\partial_{1},\partial_{2},ik)w^{1}\cdot u^R - w^{1}\cdot  {\cal N}(\partial_{1},\partial_{2},ik)u^R\big)\,dx_2\bigg|_{x_1=-0}^{x_1=+0}.
\end{multline*}
Using the relation (\ref{S1}) in the left-hand side and the relation (\ref{8}) in the right-hand side we arrive at
\begin{equation}\label{diff}
c_v^{+}=\frac {c^+_v - c_v^-}2= \frac{\int\limits_{0}^{\infty}\left(\sigma_{22}[u^R]\overline{u^R_{2}}+\sigma_{32}[u^R]\overline{u^R_{3}}-u^R_2\overline{\sigma_{11}[u^R]}\right)dx_2}
{\int\limits_{0}^{\infty}\left(\sigma_{21}[v^{R}]\overline{u^R_{2}}+\sigma_{31}[v^{R}]\overline{u^R_{3}}-v^{R}_{1}\overline{\sigma_{11}[u^R]}\right)dx_2}.
\end{equation}
Due to the explicit formulae (\ref{R1}) a direct calculation of (\ref{diff}) gives
\begin{equation}\label{diff2}
c^+_v=-\,\frac{k(1-B)^{\frac{1}{2}}(4-B)(2-B)^2(8(1-B)^2+B^2(2-B))}
{2B(8(1-B)+B^{2})(8(1-B)^2+B^2(3-2B))}.
\end{equation}
Since $B\in(0,1),$ we evidently have $c_v^{+}<0$, this inequality will be used below. Similarly, applying Green's formula to $w^1$ and $v^R$ we obtain $c_u^+ =0$. 
\medskip

The next terms of the expansion (\ref{ansatz}) can also be found; however, we do not need their explicit expressions.

\section{The outer expansions}\label{outer}

\subsection{Construction of expansions}

To construct the ``outer'' asymptotic expansion of $u^{\e}$ which is valid for large $|x_1|$, we  need to find special solutions of the problem (\ref{2dim11}) in the half-plane $\Omega^0$
\begin{equation}\label{Spec}
W^{\e\pm}(x_1,x_2)=e^{\xi^{\pm}x_1}U^{\e\pm}(x_2),\quad \xi^{\pm}=\xi^{\pm}(\e),
\end{equation}
describing waves with the exponential decay as $x_1\rightarrow\pm\infty$, that is, $\pm \xi^{\pm}<0$, with the spectral parameter
\begin{equation}\label{Lambda}
\omega^2=\omega_R^2-\e^2\Lambda^{1}+\widetilde{\Lambda}^{\e},\quad \Lambda^{1}>0,\quad \widetilde{\Lambda}^{\e}=O(\e^3).
\end{equation}

Substitution of (\ref{Spec}) into (\ref{1}) shows that $(-\xi^{\pm},U^{\e\pm})$ is a nontrivial solution of the problem
\begin{equation}\label{ODE4}
\aligned
{\cal L}(\xi^{\pm},\partial_2,ik)U^{\e\pm}(x_2)=\omega^2 U^{\e\pm}(x_2),&& x_2\in \mathbb{R}_+;\\
{\cal N}(\xi^{\pm},\partial_2,ik)U^{\e\pm}(0)=0.
\endaligned
\end{equation}

Following \cite[Ch. 9]{VT}\footnote{This book deals with linear but non-self-adjoint pencils and reduction of a quadratic pencil to a linear one is obvious.}, we accept the standard asymptotic ansatz
\begin{equation}\label{An1}
\xi^{\pm}=\mp\e \xi^{1}+\widetilde{\xi}^{\pm}, \qquad \xi^{1}>0,\qquad |\widetilde{\xi}^{\pm}|\le C\e^{2};
\end{equation}
\begin{equation}\label{An2}
U^{\e\pm}(x_2)=U^0(x_2)\mp\e \xi^{1}U^{1}(x_2)+\widetilde{U}^{\e\pm}(x_2),\qquad \|\widetilde{U}^{\e\pm}(x_2);\big(H^2(\mathbb{R}_+)\big)^3\|\le C\e^{2},
\end{equation}
with small remainders $\widetilde{\xi}^{\pm}$ and $\widetilde{U}^{\e\pm}$; their estimates are provided by general results in \cite[Ch. 9]{VT} after performing the necessary calculation as below.

We insert (\ref{Lambda}), (\ref{An1}), and (\ref{An2}) into (\ref{ODE4}) and collect coefficients at %equal 
powers of the small parameter $\e$. At the first step we obtain the problem (\ref{ODE1}) that gives $U^{0}(x_2)=u^R(x_2)$ (up to a multiplicative constant).

At the second step, cancelling the common factor $\xi^{1}$ we obtain the problem (\ref{ODE2}) and conclude that $U^{1}(x_2)=\big(-ik^{-1}u^R_{3}(x_2),0,0\big)$.

At the third step, we have (compare with (\ref{ODE3}))
\begin{equation}\label{Pr3}
\aligned
 {\cal L}(0,\partial_{2},ik)U^{2}(x_2)\, -&\, \omega^{2}_{R}U^{2}(x_2)=-\Lambda^{1}U^{0}(x_2)
\\
-&\,(\xi^{1})^2\big(\partial_{\xi}{\cal L}(0,\partial_{2},ik)U^{1}(x_2)+\tfrac{1}{2}\partial^{2}_{\xi} {\cal L}(0,\partial_{2},ik)U^0(x_2)\big),\quad x_2>0;\\
 {\cal N}(0,\partial_{2},ik)U^{2}(0)=&-(\xi^{1})^2 \big(\partial_{\xi}{\cal N}(0,\partial_{2},ik)U^{1}(0)+\tfrac{1}{2}\partial^{2}_{\xi} {\cal N}(0,\partial_{2},ik)U^0(0)\big).
\endaligned
\end{equation}
To derive the only compatibility condition, we multiply the equation in (\ref{Pr3}) by $U^0$, integrate over ${\mathbb R}_+$ and integrate by parts. These give
\begin{equation}\label{CC}
b(\xi^{1})^{2}+\Lambda^{1}\|U^{0};\big(L^2(\mathbb{R}_+)\big)^3\|^{2} =0,
\end{equation}
where $b$ is defined in (\ref{b=0}).

Since $b<0$ according to (\ref{b}) and (\ref{Req}), equation (\ref{CC}) has two roots, positive and negative, as we needed in formulae (\ref{Spec}). We have concluded with the asymptotics (\ref{An1}) 
and (\ref{An2}) together with estimates of the remainders according to \cite[Ch. 9]{VT}.

\subsection{The matching procedure}\label{matching}

We use the method of matched asymptotic expansions \cite{vD}, \cite{Il} applied in the same way as in  \cite{SN2010}, \cite{SN2011}. Since a trapped elastic mode is defined up to a multiplicative constant we regard (\ref{ansatz}) as its inner expansion 
(in a finite part of $\Pi^\e$) and take two outer expansions (for $\pm x_1\to+\infty$) in the form 
$$
C_{\pm}(\e)W^{\e\pm}(x)+\ldots,
$$ 
where $W^{\e\pm}$ are the waves introduced in (\ref{Spec}).

Since the main terms in the expansions (\ref{ansatz}) and (\ref{An2}) are equal $w^0=U^0=u^R$, we can write $C_{\pm}(\e)=1+\e C^{1}_{\pm}+\ldots$. Using (\ref{An1}) and (\ref{An2}) we derive
\begin{equation}\label{S22}
C_{\pm}(\e)W^{\e\pm}(x)=U^{0}(x_2)+\e(\mp \xi^{1}(x_1 U^{0}(x_2)+U^{1}(x_2))+C^{1}_{\pm}U^{0}(x_2))+\ldots
\end{equation}
Comparing the terms of order $\e$ in (\ref{S22}) and (\ref{ansatz}) and taking into account (\ref{S1}) we obtain
\begin{equation}\label{S24}
\xi^{1}=-c^+_{v},\quad C^{1}_{\pm}=c^{\pm}_{u}=0.
\end{equation}
Notice that the first relation in (\ref{S24}) is justified by the above-verified inequality $c^+_v < 0$. From (\ref{CC}) and (\ref{S24}) we derive the final formula for the first correction term in the 
eigenvalue expansion (\ref{Lambda}):
\begin{equation}\label{lam}
\Lambda^{1}=\frac {|b|(c^+_{v})^2} {\|U^0;\big(L^2(\mathbb{R}_+)\big)^3\|^2}>0.
\end{equation}

\section{Justification of asymptotic formulae}\label{Justify}

\subsection{Reduction to an abstract equation}\label{reduction}

Asymptotic structures found out in previous sections allow to construct a vector function 
${\bf U}^\e\in \left(H^{1}_{\sharp}(\Pi^{\e})\right)^3$ such that
$$
a_{\Pi^\e}(ik;{\bf U}^{\e},{\bf U}^{\e})<\omega_R^2 ({\bf U}^\e,{\bf U}^\e)_{\Pi^\e},
$$
where the form $a_{\Pi^\e}(ik;\cdot,\cdot)$ is defined in (\ref{tozhd})). 
By the variational principle (see, e.g., \cite[Sec. 10.2]{BS}) this implies that the (discrete) spectrum of the operator ${\cal A}^\e (ik)$ below $\omega_R^2$ is not empty and thus proves the existence of a wedge mode.

To localize the eigenvalue of ${\cal A}^\e (ik)$ and to justify its asymptotics as $\e\to0$, we employ a well-known approach based on elementary technique of spectral measure, cf. \cite{SN2010}, \cite{SN2011}.

Recall that $\left(H^{1}_{\sharp}(\Pi^{\e})\right)^3$ is the subspace of vector functions in $\left(H^{1}(\Pi^{\e})\right)^3$ satisfying condition (\ref{3a}). Denote by $\mathcal{H}^\e$ the space 
$\left(H^{1}_{\sharp}(\Pi^{\e})\right)^3$ equipped with the new scalar product
\begin{equation}\label{def1}
\langle u,v\rangle_{\e}:=a_{\Pi^\e}(ik;u,v)+ (u,v)_{\Pi^\e},
\end{equation}
and notice that the corresponding norm is equivalent to the standard norm in $\left(H^{1}_{\sharp}(\Pi^{\e})\right)^3$, see 
\cite[Ch. III, Sec. 3]{DL} and \cite{AG}. 

The quadratic form $(u,u)_{\Pi^\e}$ generates a bounded positive self-adjoint operator ${\cal K}^\e$ in $\mathcal{H}^\e$, and the problem (\ref{tozhd}) is equivalent to the abstract equation
\begin{equation}\label{def2}
{\cal K}^\e u^\e=\kappa u^\e 
\end{equation}
with the new spectral parameter $\kappa=(1+\omega^2)^{-1}$. Since the essential spectrum of ${\cal A}^\e (ik)$ coincides with the ray $[\omega_R^2,+\infty)$, %we conclude that 
the essential spectrum of ${\cal K}^\e$ is the closed segment $[0,\kappa_\dag]$ with $\kappa_\dag=(1+\omega_R^2)^{-1}$.

The spectral measure $dE(t)$ associated with ${\cal K}^\e$ (see \cite[Ch. 5]{BS}) gives rise to the family of scalar measures
$$
de_u(t):=\langle dE(t)u,u\rangle_{\e}, \qquad u\in \mathcal{H}^\e.
$$
Moreover, the following obvious formulae hold for arbitrary $u\in \mathcal{H}^\e$:
\begin{equation*}\label{spectral_int}
\|u;\mathcal{H}^\e\|^2\equiv \langle u,u\rangle_{\e}=\int\limits_{\mathbb R} de_u(t), \qquad
\| {\cal K}^\e u-\kappa u;\mathcal{H}^\e\|^2= \int\limits_{\mathbb R} (t-\kappa)^2de_u(t).
\end{equation*}

Now we assume that an interval $\Delta=(\kappa-\delta,\kappa+\delta)$ is free of the spectrum of ${\cal K}^\e$. Then we have for any $u\in \mathcal{H}^\e$
\begin{equation}\label{free_segment}
\| {\cal K}^\e u-\kappa u;\mathcal{H}^\e\|^2= \int\limits_{\mathbb R\setminus\Delta} (t-\kappa)^2de_u(t) \ge \delta^2\int\limits_{\mathbb R\setminus\Delta} de_u(t)=\delta^2\,\|u;\mathcal{H}^\e\|^2.
\end{equation}
Thus, if we find some $\delta=\delta(\e)>0$, $\kappa=\kappa(\e)>\kappa_\dag+\delta$ and ${\bf U}^\e\in \mathcal{H}^\e$ such that
\begin{equation}\label{nonempty}
\| {\cal K}^\e {\bf U}^\e-\kappa {\bf U}^\e;\mathcal{H}^\e\| <\delta\,\|{\bf U}^\e;\mathcal{H}^\e\|,
\end{equation}
then the observed contradiction between (\ref{free_segment}) and (\ref{nonempty}) shows that the interval $\Delta$ contains a spectrum point. This point is definitely an isolated eigenvalue because $\Delta$ 
does not intersect the essential spectrum of ${\cal K}^\e$ by the choice of $\kappa$ and $\delta$.

\subsection{Construction of the trial function}

\begin{figure}[ht]
	\begin{center}
		\begin{minipage}[h]{0.49\linewidth}
			\includegraphics[width=1\linewidth]{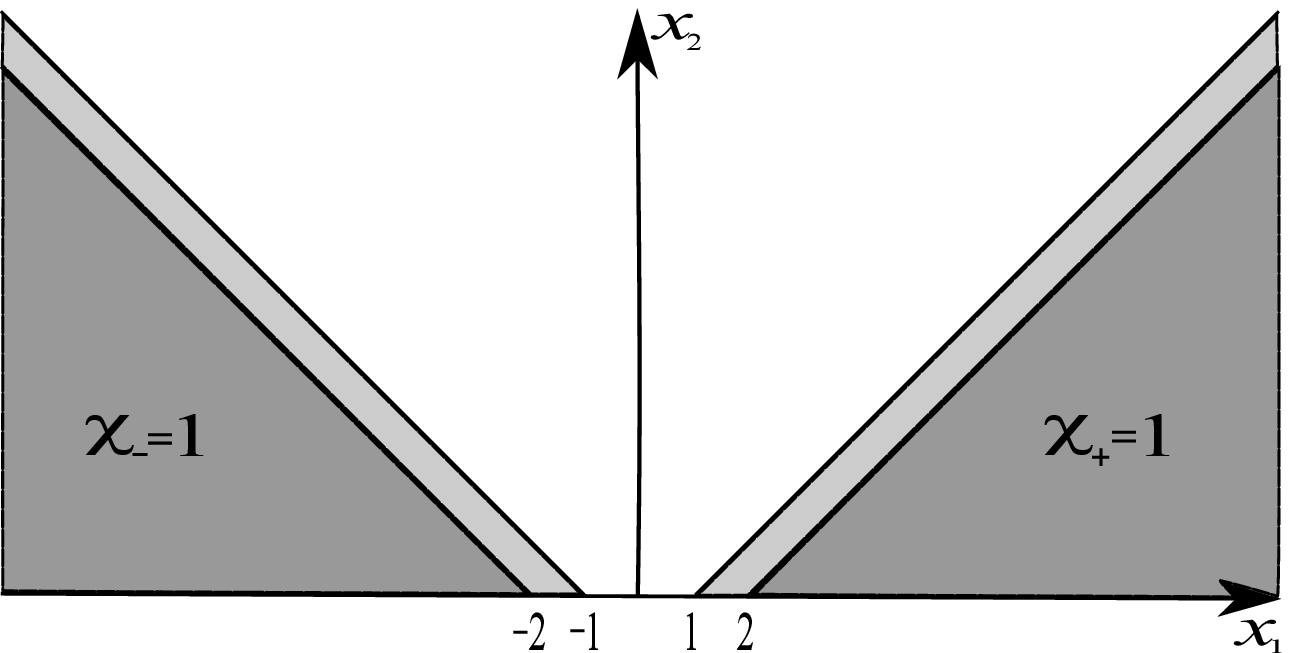}
			\caption{Cutoff functions $\chi_{\pm}$} %% 
			\label{Fig5} %% 
		\end{minipage}
		\hfill 
		\begin{minipage}[h]{0.49\linewidth}
			\includegraphics[width=1\linewidth]{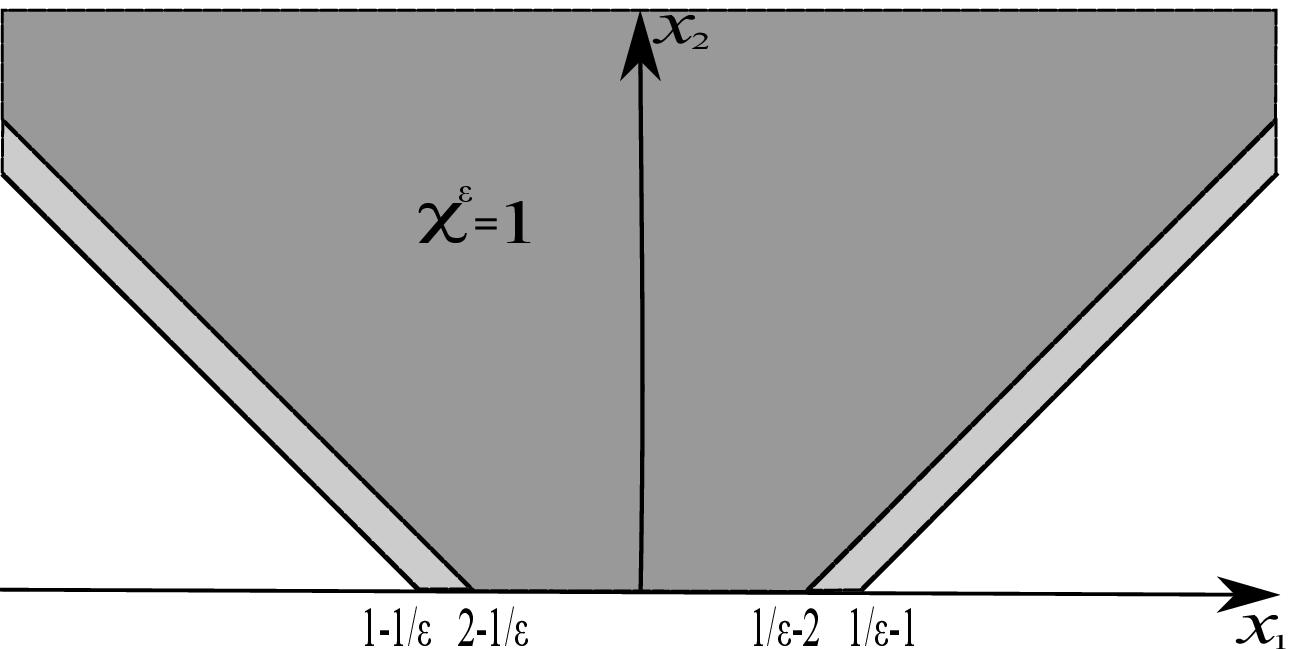}
			\caption{Cutoff function $\chi^\e$}
			\label{Fig6}
		\end{minipage}
	\end{center}%% 
\end{figure}

According to (\ref{Lambda}) and (\ref{def2}), we take 
\begin{equation}\label{spb1}
\kappa=\kappa(\e)=(1+\omega_R^2-\e^2\Lambda^1)^{-1}>\kappa_\dag+\e^2\Lambda^1\kappa_\dag^2
\end{equation}
as an approximate eigenvalue of the operator
${\cal K}^\e$, i.e. as a candidate to fulfil the inequality (\ref{nonempty}). The construction of corresponding approximate eigenvector ${\bf U}^\e$ is based on our previous asymptotic analysis and glues the 
inner and outer expansions (\ref{ansatz}), (\ref{Spec}) as follows: 
\begin{equation}\label{spb3}
{\bf U}^\e(x)={\cal X}^\e(x)W^\e(x)+\sum_{\pm} {\cal X}_\pm(x) W^{\e\pm}(x)-\sum_{\pm}{\cal X}_\pm(x){\cal X}^\e(x)(u^R(x_2)+\e c_v^{\pm}v^R(x)).
\end{equation}
Here ${\cal X}^\e$ and ${\cal X}_\pm$ are cutoff functions with overlapping supports:
\begin{equation}\label{spb2}
{\cal X}^\e(x)=\chi(x_2-|x_1|+\e^{-1}), \qquad 
{\cal X}_\pm(x)=\chi(\pm x_1-x_2),
\end{equation}
where the cutoff function $\chi$ is introduced in (\ref{chi}), see Fig. \ref{Fig5} and Fig. \ref{Fig6}.

Next, the vector function $W^\e$ has the form
$$
W^\e(x)=w^0(x)+\e w^1(x)+\e^2 w_\bullet^\e(x),
$$
where $w^0=u^R$ and $w^1$ are two terms of the inner expansion (\ref{ansatz}) constructed in Section~\ref{inner}, and $w_\bullet^\e$ is a smooth vector function supported in the union of two skewed semi-infinite strips 
$$
\Upsilon^\e_\pm=\left\{(x_1,x_2)\in \mathbb R^2:\ x_2\ge0,\ 0\le\pm x_1- \e x_2\le2\right\}
$$
and such that $W^\e$ verifies the transmission condition (\ref{3a}). Furthermore, $W^{\e\pm}$ are the waves (\ref{Spec}) with slow exponential decay as $x_1\rightarrow\pm\infty$, see (\ref{An1}). Finally, we notice 
that the terms matched in Section~\ref{matching} are duplicated in the first and second terms in (\ref{spb3}), and we eliminate this duplication by the subtrahend involving both cutoff functions (\ref{spb2}).

We emphasize that the supports of the cutoff functions drawn in Fig. \ref{Fig5} and Fig. \ref{Fig6} are angular. This becomes the key argument in the forthcoming estimates because the exponential decay of the ingredients 
of (\ref{spb3}) in the supports of $\nabla {\cal X}^\e$ and $\nabla {\cal X}_\pm$ is based on the exponential decay as either $x_2\to+\infty$ or $x_1\to\pm\infty$.

To estimate the norm on the left-hand side of (\ref{nonempty}), we use the obvious relation
\begin{equation}\label{spb5}
\|u;\mathcal{H}^\e\| =\sup\limits_{\|v;\mathcal{H}^\e\|\le1}|\langle u,v\rangle_{\e}|. 
\end{equation}
The equivalence of (\ref{tozhd}) and (\ref{def2}), definition (\ref{spb1}) and integrating by parts yield
\begin{equation}\label{spb7}
\aligned
&\langle \kappa{\bf U}^\e- {\cal K}^\e{\bf U}^\e,v\rangle_{\e}=\kappa \,\big(a_{\Pi^\e}(ik;{\bf U}^\e,v)- (\omega_R^2-\e^2\Lambda^1)({\bf U}^\e,v)_{\Pi^\e}\big)\\
=&\ \kappa \,\big(({\cal L}(\partial_1,\partial_2,ik)-\omega^2_R+\e^2\Lambda^1){\bf U}^\e,v\big)_{\Pi^\e}+
\kappa \int\limits
_{\partial\Pi^\e}{\cal N}(\partial_1,\partial_2,ik){\bf U}^\e\cdot v\,d\Gamma.
\endaligned
\end{equation}

We use for brevity the notation ${\cal L}={\cal L}(\partial_1,\partial_2,ik)$ and write
\begin{equation}\label{spb8}
\aligned
({\cal L}-\omega^2_R+\e^2\Lambda^1){\bf U}^\e= &\ {\cal X}^\e({\cal L}-\omega^2_R)(w^0+\e w^1)+
\e^2 {\cal X}^\e ({\cal L}-\omega^2_R+\e^2\Lambda^1)w_\bullet^\e\\
+ &\ \e^2{\cal X}^\e \Lambda^1\Big(w^0+\e w^1-\sum_{\pm}{\cal X}_\pm(u^R+\e c_v^{\pm}v^R)\Big)\\
+ &\, \big[{\cal L},{\cal X}^\e\big]\Big(w^0+\e w^1-\sum_{\pm}{\cal X}_\pm(u^R+\e c_v^{\pm}v^R)\Big)\\
+ &\, \sum_{\pm} {\cal X}_\pm ({\cal L}-\omega^2_R+\e^2\Lambda^1)W^{\e\pm}\\
+ &\, \sum_{\pm}\big[{\cal L},{\cal X}_\pm\big] \big(W^{\e\pm}
-{\cal X}^\e(u^R+\e c_v^{\pm}v^R)\big)\\
=: &\ J_1+J_2+J_3+J_4+\sum_{\pm}J_5^{\pm}+\sum_{\pm}J_6^{\pm},
\endaligned
\end{equation}
where $\big[{\cal L},{\cal X}\big]$ stand for the commutators of the Lame operator with cutoff functions.

We estimate the contribution of all terms to (\ref{spb7}) separately.

{\bf 1}$^\circ$. By definition of $w^0$, $w^1$ and $W^{\e\pm}$ we have $J_1\equiv0$ and $J_5^{\pm}\equiv0$.

{\bf 2}$^\circ$. Notice that the discrepancy of $w^0+\e w^1$ in the transmission conditions (\ref{3a}), (\ref{4a}) are functions of order $\e^2$ exponentially decaying in $x_2$, see (\ref{nu}), (\ref{R1}) and (\ref{8}).  By Remark \ref{Psi}, $w^0+\e w^1\in \big(W^{2,p}_{\rm{loc}}(\Pi^\e)\big)^3$ for $p<2$. 
Therefore we can fulfil the estimate 
\begin{equation}\label{spb4}
\|w_\bullet^\e;\big(W^{2,p}(\Pi^\e)\big)^3\| =
\|w_\bullet^\e;\big(W^{2,p}(\Upsilon^\e_+\cap\Upsilon^\e_-)\big)^3\| \le c
\end{equation}
(here and in the sequel we denote by $c$ constants independent of $\e$). This yields
$$
\sup\limits_{\|v;\mathcal{H}^\e\|\le1}|(J_2,v)_{\Pi^\e}| \le \e^2 \sup\limits_{\|v;\mathcal{H}^\e\|\le1} \|w_\bullet^\e;\big(W^{2,p}(\Pi^\e)\big)^3\|\,\|v;\big(L^{p'}(\Pi^\e)\big)^3\|\le c\e^2 
$$
(here $p'$ is the H\"{o}lder conjugate exponent to $p$).

{\bf 3}$^\circ$. Taking into account the supports of cutoff functions (\ref{spb2}) we have
$$
J_3=\e^2 \Lambda^1(1-{\cal X}_+-{\cal X}_-)u^R+ \e^3{\cal X}^\e \Big(w^1-\sum_{\pm}{\cal X}_\pm c_v^{\pm}v^R\Big).
$$
Using again the position of supports of cutoff functions and (\ref{S1}) we see that the both terms are of the exponential decay as $|x|\to\infty$. Thus,
$$
\sup\limits_{\|v;\mathcal{H}^\e\|\le1}|(J_3,v)_{\Pi^\e}| \le c\e^2. 
$$

{\bf 4}$^\circ$. Since the coefficients of the operator $\big[{\cal L},{\cal X}^\e\big]$ are supported in the set 
$\text{supp}(\nabla {\cal X}^\e)$ where ${\cal X}_{\pm}\equiv1$, we have
$$
J_4=\e\big[{\cal L},{\cal X}^\e\big]\widetilde{w}^{1},
$$
where $\widetilde{w}^{1}$ is the remainder term in (\ref{S1}) decaying exponentially as $|x|\to\infty$. Therefore, 
$$
\sup\limits_{\|v;\mathcal{H}^\e\|\le1}|(J_4,v)_{\Pi^\e}| \le e^{-c\e^{-1}}. 
$$

{\bf 5}$^\circ$. From relations (\ref{An1}) and (\ref{An2}), we derive that the functions $W^{\e\pm}-{\cal X}^\e(u^R+\e c_v^{\pm}v^R)$ decay exponentially as $|x|\to \infty$ on the supports of $\big[{\cal L},{\cal X}_\pm\big]$, 
respectively. Furthermore, they are of order $\e^2$, and we derive
$$
\sup\limits_{\|v;\mathcal{H}^\e\|\le1}|(J_6^{\pm},v)_{\Pi^\e}| \le c\e^2. 
$$

Next, we set ${\cal N}={\cal N}(\partial_1,\partial_2,ik)$ and rewrite the last term in (\ref{spb7}) as follows:
\begin{equation}\label{boundary}
\aligned
{\cal N}{\bf U}^\e= &\ {\cal X}^\e {\cal N}(w^0+\e w^1)+\e^2 {\cal X}^\e  {\cal N}w_\bullet^\e\\
+ &\, \big[{\cal N},{\cal X}^\e\big]\Big(w^0+\e w^1-\sum_{\pm}{\cal X}_\pm(u^R+\e c_v^{\pm}v^R)\Big)\\
+ &\, \sum_{\pm} {\cal X}_\pm {\cal N}W^{\e\pm}%\\
+ %&\, 
\sum_{\pm}\big[{\cal N},{\cal X}_\pm\big] \big(W^{\e\pm}
-{\cal X}^\e(u^R+\e c_v^{\pm}v^R)\big)\\
=: &\ J_7+J_8+J_9+\sum_{\pm}J_{10}^{\pm}+\sum_{\pm}J_{11}^{\pm}.
\endaligned
\end{equation}
We estimate the contribution of all terms separately taking into account the supports of cutoff functions (\ref{spb2}).

{\bf 6}$^\circ$. By definition of $w^0$, $w^1$ and $W^{\e\pm}$ we have $J_{10}^{\pm}\equiv0$ and
$$
\int\limits_{\partial\Pi^\e}J_7\cdot v\,d\Gamma=\int\limits_{\{\pm x_1=\e x_2>0\}}{\cal N}(w^0+\e w^1)\cdot v\,d\Gamma.
$$
We employ the standard trace inequalities
\begin{equation}\label{trace}
\aligned
& \|u;\big(H^1(\partial\Pi^\e)\big)^3\|\le c\|u;\big(W^{2,p}(\Pi^\e)\big)^3\|\quad  \text{for any}\quad   p\ge\frac 43,\\
& \|v;\big(L^2(\partial\Pi^\e)\big)^3\|\le c\|v;\mathcal{H}^\e\|,
\endaligned
\end{equation}
and recall that the discrepancy of $w^0+\e w^1$ in the transmission conditions (\ref{3a}), (\ref{4a}) are functions of order $\e^2$ exponentially decaying in $x_2$. By Remark \ref{Psi}, this gives the estimate
$$
\sup\limits_{\|v;\mathcal{H}^\e\|\le1}\Big|\int\limits_{\partial\Pi^\e}J_7\cdot v\,d\Gamma\Big|\le c\e^2.
$$

{\bf 7}$^\circ$. Similarly, from (\ref{spb4}) and (\ref{trace}) we readily have
$$
\sup\limits_{\|v;\mathcal{H}^\e\|\le1}\Big|\int\limits_{\partial\Pi^\e}J_8\cdot v\,d\Gamma\Big|\le c\e^2.
$$

{\bf 8}$^\circ$. The terms $J_9$ and $J_{11}^{\pm}$ are processed in the same way as $J_4$ and $J_6^{\pm}$, respectively, so that
$$
\sup\limits_{\|v;\mathcal{H}^\e\|\le1}\Big|\int\limits_{\partial\Pi^\e}J_9\cdot v\,d\Gamma\Big|\le e^{-c\e^{-1}};\qquad
\sup\limits_{\|v;\mathcal{H}^\e\|\le1}\Big|\int\limits_{\partial\Pi^\e}J_{11}^{\pm}\cdot v\,d\Gamma\Big|\le c\e^2.
$$
Summing up the estimates obtained, we conclude from (\ref{spb5}) and (\ref{spb7})
\begin{equation}\label{lhs}
\|\kappa{\bf U}^\e- {\cal K}^\e{\bf U}^\e;\mathcal{H}^\e\|\le c\e^2. 
\end{equation}

To estimate of the right-hand side in (\ref{nonempty}) from below, we write
\begin{equation}\label{rhs}
\|{\bf U}^\e;\mathcal{H}^\e\|^2\ge \|{\bf U}^\e;L^2(\Pi^\e)\|^2\ge \int\limits_{\e^{-1}}^{\infty} \int\limits_0^1|W^{\e+}(x)|^2\,dx_2dx_1\ge c \int\limits_{\e^{-1}}^{\infty} \! e^{2\xi^+x_1}\,dx_1\ge c\e^{-1}.
\end{equation}
Inequalities  (\ref{lhs}) and (\ref{rhs}) provide (\ref{nonempty}) with $\delta\le c\e^{\frac 52}$. According to (\ref{spb1}) we conclude that, for all sufficiently small $\e>0$, the interval 
$$
\Delta(\e)=(\kappa-c\e^{\frac 52},\kappa+c\e^{\frac 52})
$$
contains an eigenvalue $\widehat\kappa=\widehat\kappa(\e)$ of the operator ${\cal K}^\e$. Recalling the relation between $\kappa$ and $\omega$ we see that the problem (\ref{tozhd}) has an eigenvalue 
$\widehat\omega^2=\widehat\omega^2(\e)=\widehat\kappa^{-1}-1$ satisfying the asymptotic formula
\begin{equation}\label{eigen}
\widehat\omega^2=\omega_R^2-\e^2\Lambda^1+O(\e^{\frac 52}),\qquad \e\to0.
\end{equation}

\subsection{Proof of Theorem \ref{main}}

The existence of the solution $u^\e$ to the eigenproblem (\ref{2dim11}) is proved in the previous subsection. Since the Rayleigh wave, the first approximation in the asymptotic expansions (\ref{ansatz}) and (\ref{S22}), 
satisfies the symmetry condition (\ref{symm}), it is easy to observe that all other terms of expansions are also symmetric, whence $u^\e$ is symmetric. Its exponential decay follows from the relation 
$\widehat\omega^2<\omega_R^2$, see \cite{Shn} and \cite[\S53]{Gl} for scalar equations. Furthermore, by standard elliptic theory $u^\e$ is infinitely smooth outside the origin. By Remark~\ref{Psi}, $u^\e\in\big(H^2(\Omega^{\e})\big)^3$. 

Formula (\ref{asymp-velo}) with $\vartheta=\Lambda^1\omega_R^{-2}$ follows from (\ref{eigen}). Finally, we recall that $\Lambda^1$ is given by (\ref{lam}). A direct calculation using explicit formulae (\ref{R1}) gives
$$
\|U^0;\big(L^2(\mathbb{R}_+)\big)^3\|^2=\frac{B^{2}((1-B)^{2}(7-2B)+1)}{8k(1-B)^{3/2}(2-B)^2},
$$
and due to (\ref{diff2}) and (\ref{b})  $\vartheta$ depends on $\sigma$ only. 
\hfill$\square$

\begin{Remark}
Notice that despite the fact that $u^\e\in\big(H^2(\Omega^{\e})\big)^3$, the constructed approximate eigenvector ${\bf U}^\e$ lives in $\big(W^{2,p}(\Omega^{\e})\big)^3$ for $p<2$ only, see Remark \ref{Psi}. It could be smoothed out 
near the corner point with the help of an approach as in \cite{naJVMMF} but going over to the abstract equation (\ref{def2}) obviates this complication.
\end{Remark}

\section{Conclusions}

We prove the existence of a symmetric wedge mode in an elastic deformable wedge for all admissible values of the Poisson ratio $\sigma\in(-1,\frac 12)$ and interior angles close to $\pi$ and derive an 
asymptotic formula for corresponding eigenvalue. A variant of the method of matched asymptotic expansions is used. To construct the ``inner'' asymptotic expansion, we modify advisedly the use of the Kondrat'ev theory.

\bigskip
\noindent {\bf Acknowledgements.} Authors' work was partially supported by the Russian Foundation on Basic Research: joint RFBR--DFG project 20-51-12004 (A.N.), project 18-01-00325 (S.N.), project 20-01-00627 (G.Z.).

We are grateful to the referees and Associate Editor for many valuable comments and suggestions.

\medskip

\end{document}